\documentclass{amsart}
\usepackage{amsmath}
\usepackage{graphicx}
\usepackage{amsfonts}
\usepackage{amssymb}
\usepackage{hyperref}

\newtheorem{theorem}{Theorem}
\theoremstyle{plain}

\newtheorem{lemma}{Lemma}

\newtheorem{proposition}{Proposition}
\newtheorem{remark}{Remark}

\numberwithin{equation}{section}

\begin{document}
\title{Steady States for One Dimensional Conformal Metric Flows}
\author  {Yilong Ni and Meijun Zhu}
\address{Department of Mathematics\\
The University of Oklahoma\\
Norman, OK 73019\\
}
\begin{abstract}
We define two conformal structures on $S^1$ which give rise to a
different view of the affine curvature flow and a new curvature
flow, the ``$Q$-curvature flow". The steady state of these flows
are studied. More specifically,  we prove four sharp inequalities,
which state the existences of the corresponding  extremal metrics.
\end{abstract}

\maketitle

\section{Introduction}
The study of geometric flow equations have attracted great
attentions recently from scientific community. In its very short
history of existence, there have been  some very successful and
beautiful results.  To name a few, we mention  the study of
Thurston's geometrization via Ricci flow and  the applications of
geometric flows in image processing. The philosophy underneath might
be that, while mathematicians like ``stabilized'' manifolds, human vision 
system {\it prefers} ``stabilized''  images,
which could be explained as the limits of certain flow processes.

In this paper we shall start to study certain one dimensional
geometric flows, which are directly linked to image processing. Our
approach however is inspired and strongly influenced by the study
of flow equations for  high dimensional manifolds, in particular
by the study of Yamabe flow for dimension $n\ge 3$.

Recall that  on a $n$-dimensional compact Remannian manifold
$(M^n,g)$($n\ge 3$), the conformal Laplacian is defined as
\begin{equation}
\label{8-12-01}
L_g=\frac{4(n-1)}{n-2}\Delta_g+R_g,
\end{equation}
where $R_g$ is the scalar curvature of metric $g$. It is a
conformal covariant, namely for $\varphi>0$,
$$
L_{\varphi^{\frac{4}{n-2}}g}u= \varphi^{-\frac{n+2}{n-2}}L_g(\varphi u).
$$
The normalized Yamabe flow, so designed as to find a metric with constant scalar curvature
in a given conformal class, is defined as
\begin{equation}
\partial_t g =(\overline {R}_g-R_g) g,
\label{add7-27-1}
\end{equation}
where $\overline {R}_g = \int_M R_g dV_g/\int_M dV_g.$ The steady
state of (\ref{add7-27-1}) is a metric with constant scalar
curvature. See Ye \cite{Ye} and recent work by Schwetlick and
Struwe \cite{SS}.

It is obvious that (\ref{8-12-01}) is not a well defined
differential operator for $n\le 2$. We shall mimic the above to
define a general {\it $\alpha$-scalar curvature} for any given
closed convex curve. Let $(S^1,g_s)$ be the unit circle with the
induced metric $g_s=d \theta \otimes d \theta$ from $R^2$. For any
metric $g$ on $S^1$ (that is, we reparametrize the circle), we
write $g:=d \sigma \otimes d \sigma=v^{-4}g_s$  for some positive
function $v$ and define a general {\it $\alpha$- scalar curvature}
of $g$ for any positive constant $\alpha$  by
$$
R^\alpha_g=v^3(\alpha v_{\theta\theta}+v).
$$
Thus $R^\alpha_{g_s}=1.$
If we define the corresponding {\it $\alpha$- conformal Lapalace-Beltrami}
operator of $g$ by
$$
L^\alpha_g=\alpha\Delta_g+R^\alpha_g,
$$
where $\Delta_g=D_{\sigma \sigma}$, we are able to show  that
$L^\alpha_g$ is a conformal covariant, and to  establish the formula of
$R^\alpha_g$ under a conformal change of metric:
\begin{proposition}
For $\varphi>0$, if $g_2=\varphi^{-4}g_1$ then
$R^\alpha_{g_2}=\varphi^3 L^\alpha_{g_1}\varphi$, and
$$
L^\alpha_{g_2}(\psi)=\varphi^3 L^\alpha_{g_1}(\psi\varphi),
\quad\forall \psi\in C^2({\bf S^1}).
$$
\label{prop1-1}
\end{proposition}

It turns out, among all positive $\alpha$, two cases of $\alpha=1$
and $\alpha=4$ are of special interest.

We first show that the affine curvature of any given convex curve
can be viewed as a {\it $1$-scalar curvature}.
 To see this, let ${\bf
x}(\theta)$ ($\theta\in [0, 2\pi]$) be a closed convex curve in
${\bf R}^2$. One can introduce a new parameter of the affine
arc-length  $\sigma$ by requiring
$$[{\bf x}_\sigma, {\bf x}_{\sigma \sigma}]=1.
$$
In fact, if $k=k(\theta)$ is the curvature of the curve, one can
define
$$
\sigma(\theta)=\int_0^{\theta} k^{-2/3} d \theta=\int_0^{s}
k^{1/3} d s,
$$
where $s$ is the parameter of the arc length. Let $v=k^{1/3}$, we have
$$
g_1:=d\sigma \otimes d\sigma=v^{-4}d\theta \otimes d\theta=v^{-4}g_s.
$$
Then the affine curvature of ${\bf x}(\theta)$ is given by:
$$
\kappa=v^3(v_{\theta\theta}+v).
$$
which coincides with $R^1_{g_1}$ which we defined above. For simplicity
we write $\kappa_g=R^1_{g}$.

The affine flow, which is successfully used in the image processing (see \cite{ST}),  is defined as
$$
{\bf x}_t(\sigma,t)={\bf x}_{\sigma\sigma}.
$$
It is shown in \cite{ST} that the above flow is equivalent to
$$
\partial_t g=-\kappa_g g.
$$
Further, it  will be  shown in our forthcoming paper \cite{NZ2} that the affine flow is equivalent to
the following one-dimensional normalized affine flow:
\begin{equation}\label{affine}
\partial_t g=(\overline{\kappa}_g-\kappa_g) g,
\end{equation}
where $ \overline{\kappa}_g={\int \kappa_g d\sigma}/{\int
d\sigma}; $ More over, along the normalized flow
$\overline{\kappa}_g$ is always increaing.  To establish the
global existence, as well as to study the long time behaviour of
the flow, it is essential to obtain the upper bound for the total
curvature $\int \kappa_g d\sigma$ (note that under the normalized
flow, total affine length is fixed). In \cite{ST}, such an upper
bound is obtained due to  Blaschke-Santal\'o  inequality, which in
turn is proved   via the affine isoperimetric inequality as well
as  the classical Minkowski's mixed volume inequality. In this
paper, we shall directly establish such an upper bound via a pure
analytic argument. In fact, we will give an analytic proof to the
following general { Blaschke-Santal\'o} type inequality, which
immediately implies that $\overline{\kappa}_g$ is bounded above
along the normalized flow.
\begin{theorem}
(General Blaschke-Santal\'o) For $u(\theta)\in H^1(S^1)$ and
$u>0$, if $u$ satisfies
$$
\int_0^{2\pi}\frac{\cos(\theta+\gamma)}{u^3(\theta)}d\theta =0
$$
for all $\gamma$, then
$$
\int_0^{2\pi}(u^2_\theta-u^2)d\theta
\int_0^{2\pi}u^{-2}(\theta)d\theta\ge -4\pi^2,
$$
and the equality holds if and only if
$$
u(\theta)=c\sqrt{\lambda^2\cos^2(\theta-\alpha)+\lambda^{-2}\sin^2(\theta-\alpha)},
$$
for some $\lambda,c>0$ and $\alpha \in [0, 2\pi)$. \label{thm1-1}
\end{theorem}
In the consequence paper \cite{NZ2}, we will show that the global existence  of the
affine flow (\ref{affine}) follows from this inequality and the
maximum principle.

\medskip

For the case of $\alpha=4$, the {\it $4$-scalar curvature} $R^4_g$
can be viewed as the scalar curvature in an analogous
one-dimensional Yamabe flow. Let $g=d \delta\otimes d
\delta=u^{-4} g_s =u^{-4} d \theta \otimes d \theta$ for some
positive function $u \in C^1(S^1)$. We define the scalar curvature
$k_g$ of metric $g$ by
$$
k_g:=u^3(4u_{\theta \theta}+u),$$ and the corresponding  {\it conformal
Lapalace Beltrami} operator by
\begin{equation} A_g:=4\Delta_g
+k_{g},  \label{1-3-1}
\end{equation}
that is $k_g=R^4_g$ and $A_g=L^4_g$. Thus $A_g$ is a conformal
covariant by Proposition \ref{prop1-1} and
$$
k_{g_2}=\varphi^3 A_{g_1}\varphi,
$$
for  $g_2=\varphi^{-4}g_1$.
The one-dimensional  Yamabe
flow is guided by
\begin{equation}
\partial_t g=-k_g g.
\label{1-4-1}
\end{equation}
It will be shown \cite{NZ2} that this flow is equivalent to the normalized
Yamabe flow:
\begin{equation}
\partial_t g=(\overline{k}_g-k_g) g,
\label{1-4-1n}
\end{equation}
where $\overline{ k}_g= {\int _gk_g dS_g}/{\int d S_g}.$

The upper bound of $\overline{k}_g$, which is essential to the
proof of the global existence of the flow, is given by the
following inequality:
\begin{theorem}
For $u(\theta)\in H^1(S^1)$ and $u>0$,
$$
\int_0^{2\pi}(u^2_\theta-\frac14u^2)d\theta
\int_0^{2\pi}u^{-2}(\theta)d\theta\ge -\pi^2,
$$
and the equality holds if and only if
$$
u(\theta)=c\sqrt{\lambda^2\cos^2\frac{\theta-\alpha}2
+\lambda^{-2}\sin^2\frac{\theta-\alpha}2}
$$
for some $\lambda,c>0$  and $\alpha \in [0, 2\pi)$. \label{thm1-2}
\end{theorem}

Theorem \ref{thm1-2} has been proved in Ai, Chou and Wei \cite{ACW}
and Hang \cite{H}, motivated by different problems. Our proof,
which is based on a local type inequality (motivated by \cite{LZ}), is
completely different to theirs. See also Remark \ref{rem6-1}.

The analytic approach to the proofs of Theorem \ref{thm1-1} and
\ref{thm1-2} enables us to study other conformal curvature flows
on $S^1$ involving higher order derivatives. The higher
dimensional analogues include the studies of $Q$-curvature (see,
e.g., Fefferman and Graham \cite{FG}) and the extremal hermitian
metrics (see, e.g. Chen \cite{C}).  For any given $g$ on $S^1$  we
write $g=v^{-4/3}g_s$, where $g_s$ is the standard metric. The
{\it
 symmetric $Q$ curvature} of $(S^1,g)$ is defined as
\begin{equation}\label{a-q}
Q^A_g=\frac19v^{\frac53}(v_{\theta\theta\theta\theta}+10v_{\theta\theta}
+9v),
\end{equation}
and the corresponding operator $P^A$ is defined as
\begin{equation}\label{a-p}
P^A_g(f)=\frac 19\Delta^2_gf+\frac {10}9\nabla_g(\kappa_g\nabla_g
f)+Q^A_gf,\qquad f\in C^4({\bf x}),
\end{equation}
where $\kappa_g$ is the {\it 1-curvature}. Similar to Proposition
\ref{prop1-1}, we will show

\begin{proposition}\label{prop1-2}
(Conformal covariance of $P^A_g$) If $g_2
=\varphi^{-\frac43}g_1$, then
$Q^A_{g_2}=\varphi^{\frac53}P^A_{g_1}\varphi$ and
$$
P^A_{g_2}\psi=\varphi^{\frac53}P^A_{g_1}(\psi\varphi),\qquad
\forall\psi\in C^{4}(S^1).
$$
\end{proposition}

The {\it symmetric $Q$-curvature} flow is introduced as
\begin{equation}\label{a-*q}
\partial_tg=(Q^A_g-\overline{Q}^A_g)g,
\end{equation}
where $ \overline{Q}^A_g={\int Q^A_gdS_g}/{\int dS_g}$, along
which
 $\overline{Q}^A_g$ is always decreasing. In fact, we will show that
$\overline{Q}^A_g$ is bounded from below by a positive constant.
\begin{theorem}\label{thm1-3}
For $u(\theta)\in H^2(S^1)$ and $u>0$, if $u$ satisfies
$$
\int_0^{2\pi}\frac{\cos^{3}(\theta+\alpha)}{u^{5/3}(\theta)}d\theta=0
$$
for all $\alpha$, then
$$
\int_0^{2\pi}(u^2_{\theta\theta}-10u^2_\theta+9u^2)d\theta
\left(\int_0^{2\pi}u^{-2/3}(\theta)d\theta\right)^3\ge C_0>0.
$$
More over, if $u_0$ is an extremal function, then the symmetric
$Q$-curvature of $u_0^{-4/3} g_s$ is a constant.
\end{theorem}
The above  might be needed to obtain the global existence and
convergence of the symmetric $Q$-curvature flow. We expect such
flow will have interesting applications in image processing
(especially in image recognition), see Remark \ref{rem6-2} after
the proof.

\medskip

At this stage, it is nature to introduce  the one-dimensional {\it
$Q$ curvature} of $(S^1, g)$ given by
\begin{equation}\label{q}
Q_g=\frac19v^{\frac53}(16v_{\theta\theta\theta\theta}+40v_{\theta\theta} +9v),
\end{equation}
and the corresponding operator $P_g$ on $({S^1},g)$ given  by
\begin{equation}\label{p}
P_g(f)=\frac {16}9\Delta^2_gf+\frac{40}9\nabla_g(k_g\nabla_g
f)+ Q_gf,\qquad f\in C^4({S^1}),
\end{equation}
where $k_g$ is the {\it 4-scalar curvature} of $({S^1},g)$. The
$Q$-curvature flow is then guided by
\begin{equation}\label{*q}
\partial_tg=(Q_g-\overline{Q}_g)g,
\end{equation}
where $ \overline{Q}_g={\int Q_gdS_g}/{\int dS_g}. $ Simliar to
Proposition \ref{prop1-2}, we will show that
\begin{proposition}(Conformal Covariant of $P_g$)
If $g_2=\varphi^{-\frac43}g_1$, then
$Q_{g_2}=\varphi^{\frac53}P_{g_1}\varphi$ and
$$
P_{g_2}\psi=\varphi^{\frac53}P_{g_1}(\psi\varphi),\qquad
\forall\psi\in C^{4}({\bf x}).
$$
\label{prop1-3}
\end{proposition}

\begin{remark}
We may define general $\alpha - Q$ curvature on $(S^1,g)$ as
$$
Q^\alpha_g=v^{\frac53}(\frac{\alpha^2}{9}v_{\theta\theta\theta\theta}
+\frac{10\alpha}{9}v_{\theta\theta}+v),
$$
and the corresponding operator as
$$
P_g^\alpha(f)=\frac{\alpha^2}{9}\Delta^2_gf+\frac{10\alpha}{9}
\nabla_g(R^\alpha_g\nabla_g f)+Q^\alpha_g f,
$$
where $g=v^{-4/3}g_s$ and $R^\alpha_g$ is the $\alpha-$scalar curvature of
$(S^1,g)$. Similar to Proposition \ref{prop1-3} we have
\begin{proposition}(Conformal Covariant of $P_g^\alpha$)
If $g_2=\varphi^{-\frac43}g_1$, then
$Q^\alpha_{g_2}=\varphi^{\frac53}P^\alpha_{g_1}\varphi$ and
$$
P^\alpha_{g_2}\psi=\varphi^{\frac53}P^\alpha_{g_1}(\psi\varphi),\qquad
\forall\psi\in C^{4}({\bf x}).
$$
\label{prop1-3p}
\end{proposition}
\label{rem1-2}
\end{remark}
The following inequality, which classifies the steady state of the
flow, is essential to the proof of the global existence and
convergence of the flow.
\begin{theorem}
For $u(\theta)\in H^2(S^1)$ and $u>0$,
$$
\int_0^{2\pi}(u^2_{\theta\theta}-\frac52u^2_\theta
+\frac9{16}u^2)d\theta
\left(\int_0^{2\pi}u^{-2/3}(\theta)d\theta\right)^3\ge 9\pi^4,
$$
and the equality holds if and only if
$$
u(\theta)=c\left(\lambda^2\cos^2\frac{\theta-\alpha}2
+\lambda^{-2}\sin^2\frac{\theta-\alpha}2\right)^{3/2},
$$
for some $\lambda,c>0$ and $\alpha \in [0, 2\pi)$. \label{thm1-4}
\end{theorem}

The main theme of this paper is to prove the above four theorems.
Note that extremal functions of these four inequalities give the
extremal metrics for corresponding conformal curvatures. Therefore
we also have  classifications of the limiting metrics of the
affine flow, Yamabe flow and the $Q$-curvature flow. The steady
state for symmetric $Q$-curvature flow is not completely clear
yet.
 See Remark \ref{rem6-2}.

During our writing of this paper, we learned that the analytic
approach to the affine isoperimetric as well as the
Blasche-Santal\'o inequalities has been proposed and tried by
another group of researchers, see \cite{CHLYZ}. For example, the
analytic proof of Theorem \ref{thm1-1} was given in \cite{CHLYZ}.
However, their proof is incomplete (Lemma 6.2 in \cite{CHLYZ} is
not correct since $u_0$ can vanish on a whole interval). 
When the paper was in circulation, we were informed by J. Wei \cite{wei}
that Proposition A.5 in \cite{ACW} also yields an analytical proof of Theorem 1.
The proofs of Theorem \ref{thm1-2} and \ref{thm1-4} are also obtained
in a recent preprint \cite{H} by Hang. His work
seems to be inspired by early work of Yang and Zhu \cite{YZ} on
three dimensional Paneitz operator. Our proofs are different to
those in \cite{YZ} and \cite{H}.

 The paper is organized as follows.
In section 2, we prove the conformal covariance of operators
$L^\alpha_g$, $P^A_g$ and $P_g$. In section 3, we establish a
local sharp inequality, which implies Theorem \ref{thm1-1} using the
stereographic projection. In section 4, we give a new proof of the
Blaschke-Santal$\acute{\mbox{o}}$ inequality. The idea is also
used in the proof of Theorem \ref{thm1-3}. We prove Theorem \ref{thm1-4}
and \ref{thm1-3} in section 5.

\section{Conformal Covariance}

Suppose that $g_1:=d \sigma_1 \otimes d \sigma_1=v^{-4}g_s:=v^{-4}
d \theta \otimes d \theta$, where $g_s$ is the standard metric on
the unit circle. We will use convention that
$\Delta_{g_i}=D_{\sigma_i \sigma_i}$ and
$\nabla_{g_i}=D_{\sigma_i}$ for $i=1$ and $2$.  Then for
$g_2=\varphi^{-4} g_1=(\varphi v)^{-4} g_s$, by definition
\begin{align*}
R^\alpha_{g_2}= &(\varphi v)^3(\alpha\Delta_{g_s}(\varphi v)+ \varphi v)\\
= &(\varphi
v)^3(\alpha v\Delta_{g_s}\varphi+2\alpha\nabla_{g_s}\varphi\nabla_{g_s}v
+\alpha\varphi\Delta_{g_s}v+\varphi v).
\end{align*}
On the other hand,
\begin{align*}
\varphi^3L^\alpha_{g_1}\varphi=&\varphi^3(\alpha\Delta_{g_1}\varphi
+v^3(\alpha\Delta_{g_s}v+v)\varphi)\\
=&\varphi^3(\alpha v^2\nabla_{g_s}(v^2\nabla_{g_s}\varphi)+\alpha v^3
\varphi\Delta_{g_s}v+v^4\varphi)\\
=&(\varphi v)^3(\alpha v\Delta_{g_s}\varphi
+2\alpha\nabla_{g_s}\varphi\nabla_{g_s}v
+\alpha\varphi\Delta_{g_s}v+\varphi v).
\end{align*}
Therefore $R^\alpha_{g_2}=\varphi^3L^\alpha_{g_1}\varphi$. It
follows that, for any $\psi\in C^2({S^1}),$
\begin{align*}
L^\alpha_{g_2}\psi=&\alpha\Delta_{g_2}\psi+R^\alpha_{g_2}\psi
=\alpha\varphi^2\nabla_{g_1}(\varphi^2\nabla_{g_1}\psi)
+\varphi^3\psi L^\alpha_{g_1}\varphi\\
=&\alpha\varphi^4\Delta_{g_1}\psi+2\alpha\varphi^3\nabla_{g_1}
\varphi\nabla_{g_1}\psi+\alpha\varphi^3\psi\Delta_{g_1}\varphi+
R^\alpha_{g_1}\varphi^4\psi\\
=&\varphi^3L^\alpha_{g_1}(\psi\varphi).
\end{align*}
This completes the proof of Proposition \ref{prop1-1}.

To prove Proposition \ref{prop1-2}, we recall that for
$g=v^{-4/3}g_s$, the {\it 1-scalar curvature} $\kappa_g$ has the
following form:
\begin{equation}\label{a-newk}
\kappa_g=(v^{\frac13})^3((v^{\frac13})_{\theta\theta}
+v^{\frac13})=\frac13v^{\frac13}v_{\theta\theta}-\frac29
v^{-\frac23}(v_\theta)^2+ v^{\frac43}.
\end{equation}
For any $f\in C^4({S^1})$ and $g_1=v^{-4/3}g_s$,
\begin{align*}
\nabla_{g_1}f=&f_\theta v^{\frac23},\\
\Delta_{g_1}f=&f_{\theta\theta}v^{\frac 43}+\frac23v^{\frac13}f_\theta v_\theta,\\
\Delta_{g_1}^2f=&\frac{10}9v^{\frac23}(2v^2_\theta f_{\theta\theta}
+v_\theta v_{\theta\theta}f_\theta)\\
&\, +\frac13v^{\frac53}(3vf_{\theta\theta\theta\theta}+12v_\theta
f_{\theta\theta\theta}+8v_{\theta\theta}f_{\theta\theta}
+2v_{\theta\theta\theta}f_\theta),\\
\nabla_{g_1}(\kappa_{g_1}\nabla_{g_1}f)=&\nabla_{g_1}(f_\theta
(\frac13vv_{\theta\theta}-\frac29v^2_\theta+v^2))\\
=&v^{\frac23}(v^2f_{\theta\theta}-
\frac29v^2_\theta f_{\theta\theta}-\frac19v_\theta v_{\theta\theta}
f_\theta)\\
&\,+\frac13v^{\frac53}(v_{\theta\theta}f_{\theta\theta} +6v_\theta
f_\theta+f_\theta v_{\theta\theta\theta}).
\end{align*}
Therefore
\begin{align*}
9\varphi^{\frac53}P^A_{g_1}(\varphi)=&\varphi^{\frac53}(
\Delta_{g_1}^2\varphi+10\nabla_{g_1}(\kappa_{g_1}
\nabla_{g_1}\varphi)+9Q^A_{g_1}\varphi)\\
=&(v\varphi)^{\frac53}(v\varphi_{\theta\theta\theta\theta}
+4v_\theta\varphi_{\theta\theta\theta}+6v_{\theta\theta}
\varphi_{\theta\theta}+4v_{\theta\theta\theta}\varphi_\theta\\
&\,+20v_\theta\varphi_\theta+10v\varphi_{\theta\theta})
+9\varphi^{\frac 8 3} Q^A_{g_1}\\
=&(v\varphi)^{\frac53}((v\varphi)_{\theta\theta\theta\theta}
+10(v\varphi)_{\theta\theta}+9v\varphi)\\
=&9Q^A_{g_2}.
\end{align*}
For any $\psi>0$, let $g_3=\psi^{-4/3}g_2$. The above calculation
shows that
$$
Q^A_{g_3}=\psi^{\frac53}P^A_{g_2}\psi.
$$
On the other hand, $g_3=\psi^{-4/3}g_2=(\varphi\psi)^{-4/3}g_1$,
thus
$$
Q^A_{g_3}=(\varphi\psi)^{\frac53}P^A_{g_1}(\varphi\psi).
$$
Comparing the above two equations, we obtain that
\begin{equation}\label{a-qc}
P^A_{g_2}\psi=\varphi^{\frac53}P^A_{g_1}(\varphi\psi)
\end{equation}
for any $\psi>0$. Finally, for general $\psi$, choose $C$
sufficiently large, such that $C+\psi>0$. Applying (\ref{a-qc}) to
$C+\psi$ and using the linearity of $P^A_g$, we obtain that
(\ref{a-qc}) holds for all $\psi$. This completes the proof of
Proposition \ref{prop1-2}. The proof of Proposition \ref{prop1-3}
is very similar. We shall skip the details.

\section{A Local Sharp Inequality}

In this section, we shall follow \cite{LZ} to establish a local
sharp inequality, which will yield the  global inequality Theorem
\ref{thm1-2}.

 For any three positive constants $a,b,r>0 $, We define
$$
D_{a,b}(-r,r)=\{f(y) \ : \ f(y)-b\in H^1_0(-r,r), \  f>0, \
 \int_{-r}^r f^{-2}(y)dy=a\}
$$
\begin{proposition}
\label{prop2-1} If $a>2rb^{-2},$ then
$$
\inf_{w\in D_{a,b}(-r,r)}\int_{-r}^r |\nabla w|^2dy = 2\tau^{1/2}
\left(\frac r\lambda-\arctan\left(\frac r\lambda\right)\right),
$$
where $\tau>0$ and $\lambda>0$ are determined by
$$
2\tau^{-1/2}\arctan(r\lambda^{-1})=a,\qquad
\frac{\lambda^2+r^2}{\lambda}\tau^{1/2}=b^2.
$$
The infimum is achieved if and only if by
$$
w(y)=\tau^{1/4}\sqrt{\frac{\lambda^2+y^2}{\lambda}}.
$$
If  $a<2rb^{-2}$, then
$$
\inf_{w\in D_{a,b}(-r,r)}\int_{-r}^r |\nabla w|^2dy= 2\tau^{1/2}
\left(-\frac r\lambda+\operatorname{arctanh}\left(\frac
r\lambda\right)\right),
$$
where $\tau<0$ and $\lambda>r$ are determined by
$$
2|\tau|^{-1/2}\operatorname{arctanh}(r\lambda^{-1})=a,\qquad
\frac{\lambda^2-r^2}{\lambda}|\tau|^{1/2}=b^2.
$$
The infimum is achieved if an only if by
$$
w(y)=|\tau|^{1/4}\sqrt{\frac{\lambda^2-y^2}{\lambda}}.
$$
If $a=2rb^{-2}$, then
$$
\inf_{w\in D_{a,b}(-r,r)}\int_{-r}^r |\nabla w|^2dy= 0.
$$
The infimum is achieved if and only if by $w=b$.
\end{proposition}

\begin{proof}
It is easy to see that: if $a\ge 2 r b^{-2}$, then $$\inf_{w \in
D_{a,b}(-r, r)} \int_{-r}^r|\nabla w|^2 d y= \inf_{w \in
D_{a,b}(-r, r), w\ge b} \int_{-r}^r|\nabla w|^2 d y;$$
 And if $a\le  2 r b^{-2}$,
$$
\inf_{w \in D_{a,b}(-r, r)} \int_{-r}^r|\nabla w|^2 d y= \inf_{w
\in D_{a,b}(-r, r), w\le b} \int_{-r}^r|\nabla w|^2 d y.$$ So, if
$\{w_i\}\subset D_{a,b}(-r,r)$ is a minimizing sequence
satisfying:
$$
\lim_{i\to\infty}\int_{-r}^r|\nabla w_i|^2dy =\inf_{w\in
D_{a,b}(-r,r)}\int_{-r}^r |\nabla w|^2dy,
$$
we may assume without loss of generality (after rearrangement)
that $w_i(-y)=w_i(y)$, and when $a>2rb^{-2}$, $w'_i(y)\le 0$ for
$y>0$, $w'_i(y)\ge 0$ for $y<0$; When $a<2rb^{-2}$, $w'_i(y)\ge 0$
for $y>0$, $w'_i(y)\le 0$ for $y<0$.

Since $\{w_i\}$ is a minimizing sequence, $\int_{-r}^r|\nabla
w_i|^2dy\le C$ for all $i$. It follows from $w_i-b \in H_0^1$ that
$||w_i||_{H^1}\le C$. Therefore, up to a subsequence,
$$
w_i\rightharpoonup w_0\mbox{ weakly in } H^1(-r,r).
$$
Using Sobolev embedding theorem we obtain that
$||w_i||_{C^{0,\frac12}(-r,r)}\le C$, and $w_i\to w_0$ in
$C^{0,\alpha}(-r, r)$ for any $\alpha\in(0,\frac12)$.

Fatou's Lemma implies
$$
\int_{-r}^rw_0^{-2}(y)dy\le
\lim_{i\to\infty}\int_{-r}^r w_i^{-2}(y)dy=a.
$$
It follows from $w_0\in C^{0,\frac12}(-r,r)$ that $w_0(y)>0$ for
any $y\in (-r,r)$. Therefore $w_0\in D_{a,b}(-r,r)$ is a
minimizer. It is easy to see that $w_0(y)=w_0(-y)$,
$w_0(r)=w_0(-r)=b$, $w_0\ge b$ if $a>2b^{-2}$ and $w_0\le b$ if
$a<2b^{-2}$. Moreover $w_0$ satisfies the Euler-Lagrange equation
\begin{equation}\label{el1}
w''(y)=\tau w^{-3}.
\end{equation}
for some constant  $\tau$.

If  $a>2rb^{-2}$, it is easy to see that $\tau>0$. Let
\begin{equation}\label{8-13-01}
v_1(y)=\tau^{\frac14}
\left(\frac{\lambda^2+y^2}{\lambda}\right)^{\frac12}.
\end{equation}
We can  verify that $v_1(y)$ satisfies (\ref{el1}), and
$v'_1(0)=0=w'_0(0)$. Choosing $\lambda=\tau^{-1/2}w^2_0(0)$, we
have $v(0)=w_0(0)$. It follows from the uniqueness of the solution
of (\ref{el1}) that $w_0(y)=v_1(y)$. The boundary condition yields
$$
\tau^{\frac14}
\left(\frac{\lambda^2+r^2}{\lambda}\right)^{\frac12}=b,
$$
and $w_0\in D_{a,b}(-r, r)$ implies
$$
a=\int_{-r}^rw_0^{-2}dy=2\tau^{-1/2}\arctan(r\lambda^{-1}).
$$
Thus the infimum is given by
$$
\int_{-r}^r|\nabla w_0|^2dy=2\tau^{1/2}\left(\frac r\lambda
-\arctan(r\lambda^{-1})\right).
$$
On the other hand, if $w_0\in D_{a, b}(-r, r)$ is a minimizer,
then $w_0 \ge b$ and satisfies (\ref{el1}). From the uniqueness of
the solution to  the ODE, we know $w_0=v_0$ with
$\lambda=\tau^{-1/2}w^2_0(0)$.

If  $a<2rb^{-2}$, then  $\tau<0$. Let
$$
v_2(y)=|\tau|^{\frac14}
\left(\frac{\lambda^2-y^2}{\lambda}\right)^{\frac12}.
$$
Again $v_2(y)$ satisfies (\ref{el1}), and $v'_2(0)=0=w'_0(0)$. We
choose $\lambda=|\tau|^{-1/2}w_0^2(0)$, so that $v_2(0)=w_0(0)$.
It follows from the uniqueness of the solution of (\ref{el1}) that
$w_0(y)=v_2(y)$. The boundary condition gives
$$
|\tau|^{\frac14}
\left(\frac{\lambda^2-r^2}{\lambda}\right)^{\frac12}=b,
$$
and $w_0\in D_{a,b}(-r, r)$ implies
$$
a=\int_{-r}^rw_0^{-2}dy=2|\tau|^{-1/2}\mbox{arctanh}( r
\lambda^{-1}).
$$
So the infimum is given by
$$
\int_{-r}^r|\nabla w_0|^2dy=2|\tau|^{1/2}\left(-\frac r\lambda
+\mbox{arctanh}(r\lambda^{-1})\right).
$$
Similarly, one can see that such $v_2$ is the only minimizer.

The case of $a=2rb^{-2}$ is trivial.

\end{proof}

\section{Proof of Theorem \ref{thm1-2}}

Our proof of Theorem \ref{thm1-2} is in the same spirit as that of
Theorem C in \cite {LZ}. For $u(\theta)\in H^1(S^1)$($\theta\in
[-\pi,\pi]$), $u>0$, without loss of generality (due to
rearrangement), we may assume that $u(\theta)=u(-\theta)$ and
$u(\theta)$ is decreasing on $[0,\pi]$.

For $\epsilon\in(0,\pi)$, let
$S^1_\epsilon=\{\theta:\theta\in[-\pi+\epsilon,\pi-\epsilon]\}$,
$$
r_\epsilon=\tan\frac{\pi-\epsilon}{2}>0, \mbox{ and }
a_\epsilon=\int_{-\pi+\epsilon}^{\pi-\epsilon}u^{-2}(\theta)
d\theta.
$$
It is easy to see that $r_\epsilon\to +\infty$ as $\epsilon\to 0$,
and $\Phi(S^1_\epsilon)=[-r_\epsilon,r_\epsilon]$, where $\Phi$ is
the stereographic projection $\Phi:S^1\to {\mathbf R}^1$, defined
as
$$
y=\Phi(\theta)=\tan\frac\theta 2 \quad \mbox{for} \quad \theta\in [-\pi,\pi].
$$
Let $w(y)=u(\Phi^{-1}(y))\varphi(y)$, where
$\varphi(y)=\sqrt{(1+y^2)/2}$. Then we have (note that
$\theta=2\arctan y$)
\begin{align*}
\int_{-r_\epsilon}^{r_\epsilon}w^{-2}dy=&
\int_{-r_\epsilon}^{r_\epsilon}u^{-2}(2\arctan y)\frac2{1+y^2}dy\\
=&\int_{-\pi+\epsilon}^{\pi-\epsilon}u^{-2}(\theta)d\theta=a_\epsilon.
\end{align*}
Therefore $w(y)\in D_{a_\epsilon, b_\epsilon}
(-r_\epsilon,r_\epsilon)$, where $b_\epsilon=w(r_\epsilon)
=u(\pi-\epsilon)\sqrt{(1+r_\epsilon^2)/2}$. Using Proposition
\ref{prop2-1}, we have
\begin{equation}\label{uselemma}
\int_{-r_\epsilon}^{r_\epsilon}(w'(y))^2dy\ge 2\sqrt{\tau}
\left(\frac{r_\epsilon}{\lambda}-\arctan\frac {r_\epsilon}\lambda)
\right),
\end{equation}
where $\tau$ and $\lambda$ are determined by
\begin{equation}
\label{ab}
2\tau^{-1/2}\arctan(r_\epsilon\lambda^{-1})=a_\epsilon,\qquad
\frac{\lambda^2+r_\epsilon^2}{\lambda}\tau^{1/2}=b_\epsilon^2.
\end{equation}
On the other hand,
$$
w'(y)=\varphi\frac{d}{dy}u(2\arctan y)+u(2\arctan y)\varphi'(y)
=u'(\theta)\frac1{\varphi(y)}+u(\theta)\varphi'(y),
$$
thus
\begin{align*}
\int_{-r_\epsilon}^{r_\epsilon}(w'(y))^2dy=&
\int_{-r_\epsilon}^{r_\epsilon}[\frac{u_\theta^2}{\varphi^2(y)}
+\frac{(u^2)_\theta\varphi'(y)}{\varphi(y)}+u^2\cdot (\varphi'(y))^2]dy\\
=&\int_{-\pi+\epsilon}^{\pi-\epsilon}u^2_\theta d\theta
-\int_{-r_\epsilon}^{r_\epsilon}\varphi\varphi''u^2dy
+\varphi\varphi'u^2(2\arctan y)|_{-r_\epsilon}^{r_\epsilon}\\
=&\int_{-\pi+\epsilon}^{\pi-\epsilon}(u^2_\theta-\frac14u^2)d\theta
+b_\epsilon^2\frac{2r_\epsilon}{1+r_\epsilon^2}.
\end{align*}
Therefore (\ref{uselemma}) and (\ref{ab}) imply
\begin{align*}
\int_{-\pi+\epsilon}^{\pi-\epsilon}(u^2_\theta-\frac14u^2)d\theta
\int_{-\pi+\epsilon}^{\pi-\epsilon}u^{-2}d\theta
=&a_\epsilon(\int_{-r_\epsilon}^{r_\epsilon}(w'(y))^2dy
-b_\epsilon^2\frac{2r_\epsilon}{1+r_\epsilon^2})\\
\ge& a_\epsilon\left(2\sqrt{\tau}
\left(\frac{r_\epsilon}{\lambda}-\arctan\frac {r_\epsilon}\lambda
\right)-b_\epsilon^2\frac{2r_\epsilon}{1+r_\epsilon^2}\right)\\
=&4\arctan\frac{r_\epsilon}{\lambda}\left(
\frac{r_\epsilon}{\lambda}-\arctan\frac{r_\epsilon}{\lambda}
-\frac{r_\epsilon}{\lambda}\cdot
\frac{\lambda^2+r^2_\epsilon}{1+r_\epsilon^2}
\right)\\
=&4\arctan\frac{r_\epsilon}{\lambda}\left(
\frac{r_\epsilon}{\lambda}\cdot \frac{1-\lambda^2}{1+r_\epsilon^2}
-\arctan\frac{r_\epsilon}{\lambda}\right).
\end{align*}

Letting $\beta_\epsilon=\arctan\frac {r_\epsilon}\lambda\in(0,\frac\pi2)$,
and using (\ref{ab}), we have
$$
a_\epsilon b_\epsilon^2=2\arctan\frac{r_\epsilon}{\lambda} \cdot
\frac{\lambda^2+r^2_\epsilon}{\lambda}=
4r_\epsilon\frac{\beta_\epsilon}{\sin(2\beta_\epsilon)}.
$$
Since $b_\epsilon^2=u^2(\pi-\epsilon)(1+r^2_\epsilon)/2=O(r^2_\epsilon)$,
and $\lim_{\epsilon\to 0}a_\epsilon=\int_{-\pi}^{\pi}u^{-2}d\theta:=a$,
we have
$
\lim_{\epsilon\to 0}\beta_\epsilon=\pi/2.$
It then follows from (\ref{ab}) that $
\lim_{\epsilon\to 0}\lambda={2\pi}{a^{-1}u^{-2}(\pi)}.
$
Hence
$$
\lim_{\epsilon\to 0}4\arctan\frac{r_\epsilon}{\lambda}\left(
\frac{r_\epsilon}{\lambda}\cdot\frac{1-\lambda^2}{1+r_\epsilon^2}
-\arctan\frac{r_\epsilon}{\lambda}\right)
=4\frac\pi2(0-\frac\pi2)=-\pi^2,
$$
and
$$
\int_{-\pi}^{\pi}(u^2_\theta-\frac14u^2)d\theta
\int_{-\pi}^{\pi}u^{-2}d\theta\ge -\pi^2.
$$

The equality holds if and only if $u$ satisfies the Euler-Lagrange
equation
$$
u_{\theta \theta} +\frac 1 4 u= \tau  u^{-3}
$$
for some positive constant $\tau$.
Thus $w(y)=u(\Phi^{-1}(y)) \varphi(y)>0$ satisfies
\begin{equation}
w''=\tau w^{-3}, w(y)=O(|y|), \ \mbox{as} \ |y|\to \infty.\label{8-9-01}
\end{equation}
The first integral of equation (\ref{8-9-01}) is
$$
(w')^2=c_1-\tau w^{-2}.
$$
Since $w(y)=O(|y|)$, as $|y|\to \infty$, we may assume that $w(y)$ attains
its minimum at some point $a$. It is easy to see that $w(y)$ is increasing
when $y>a$, decreasing when $y<a$, and all solutions to (\ref{8-9-01}) are
given by
$w(y)=\sqrt{\frac{\tau+c_1^2(y-a)^2}{c_1}}$, that is
$$
u(\theta)=c\sqrt{\lambda^2\cos^2\frac{\theta-\alpha}2
+\lambda^{-2}\sin^2\frac{\theta-\alpha}2},
$$
where $c=(4\tau)^{\frac14}$ and $\lambda=(\tau)^{\frac14}/\sqrt{c_1}$.
This completes the proof of Theorem \ref{thm1-2}.

\section{A new proof of Blaschke-Santal\'o inequality}
The original proof of Blaschke-Santal\'o inequality and affine
isoperimetric inequality involved Minkowski mixed volume
inequality for convex body, see for example \cite{Lu}.  Here we
shall give a new analytic  proof. The analytic approach enables us
to extend such inequality to the one involving higher order
derivatives (see our Theorem \ref{thm1-3}), as well as allows us
to remove the geometric constrains (i.e. the convexity of the
curve). Such an approach  was used also  in \cite{CHLYZ} aiming to
generalize the affine isoperimetric inequality. However, as we
pointed out early that their proof is incomplete.

For any given $\lambda>0$ and $\theta \in [0, 2\pi]$, let
$$
\psi_\lambda(\theta)=\sqrt{\lambda^2\cos^2\theta
+\lambda^{-2}\sin^2\theta},$$
$$ \sigma_\lambda(\theta)=\int_0^\theta
\psi_\lambda^{-2}d\theta= \left\{
\begin{array}{ll}
\arctan(\lambda^{-2}\tan\theta), & \theta \in [0, \pi/2]\\
\arctan(\lambda^{-2}\tan\theta)+\pi, & \theta \in (\pi/2, 3\pi/2]\\
\arctan(\lambda^{-2}\tan\theta)+2 \pi, & \theta \in (3\pi/2, 2 \pi)
\end{array}
\right.
$$
and define
$$
(T_\lambda u)(\theta):=u(\sigma_\lambda(\theta))\psi_\lambda(\theta).
$$
From conformal invariant properties, we know that
\begin{equation}
\int_0^{2\pi}(u^2_\theta-u^2)d\theta=
\int_0^{2\pi}((T_\lambda u)^2_\theta-(T_\lambda u)^2)d\theta,
\label{8-9-02}
\end{equation}
and
\begin{equation}
\int_0^{2\pi}u^{-2}d\theta=\int_0^{2\pi}(T_\lambda u)^{-2}d\theta.
\label{8-9-03}
\end{equation}
Also one can easily check that for all $\alpha\in[0,2\pi)$
\begin{equation}
\int_0^{2\pi}\frac{\cos(\theta+\alpha)}{(T_\lambda u)^3}d\theta
=\sqrt{\lambda^{-2}\cos^2\alpha+\lambda^2\sin^2\alpha}\cdot
\int_0^{2\pi}\frac{\cos(\theta+\tilde{\alpha})}{u^3}d\theta,
\label{8-9-04}
\end{equation}
where $\tilde{\alpha}=\sigma^{-1}_\lambda(\alpha)$.

For any $u\in H^1(S^1)$, $u>0$, let
$$
J(u):=\int_0^{2\pi}(u^2_\theta-u^2)d\theta\cdot\int_0^{2\pi}
u^{-2}d\theta.
$$
If we denote
$$
H^1_s(S^1):=\{ u \in H^1(S^1) \ : \ u>0, \ \int_0^{2\pi} \cos
(\theta+\gamma) u^{-3}(\theta)d\theta =0,\mbox{ for all }\gamma
\in[0,2\pi)\},
$$
it is obvious that
$$
\inf_{H^1_s}J(u) \le J(1)=-4 \pi^2<0.
$$
We first shall prove that $\inf_{H^1_s}J(u)$ is bounded below and  the  minimizer of $J(u)$ exists.

Suppose $\{u_i\}$ is a sequence of minimizing functions for $\inf_{H^1_s}J(u)$ with
$\max_\theta u_i(\theta)=1.$
\begin{lemma}\label{eq}
For each $f\in C^0(S^1)$, there exists a constant $a$ such that $f(a)=f(a+\pi)$.
\end{lemma}
\begin{proof}
Let $g(\theta)=f(\theta+\pi)-f(\theta)$. Then
$g(\theta)=-g(\theta+\pi)$. Intermediate value theorem immediately
yields the result.
\end{proof}
Due to Lemma \ref{eq}, without loss of generality (after
rotating), we may assume that $u_i(0)=u_i(\pi)$ and $u_i(\pi/2)\ge
u_i(3\pi/2)$. Let
$$
\lambda_i=\sqrt{\frac{u_i(\frac\pi2)}{u_i(0)}},
\mbox{ and }
w_i(\theta)=\frac{T_{\lambda_i}u_i}{||T_{\lambda_i}u_i
||_{L^\infty}}.
$$
It is easy to see from (\ref{8-9-02}) to (\ref{8-9-04}) that $\{w_i\}$ is also a minimizing sequence and
$$
w_i(0)=w_i(\frac\pi2)=\frac{\sqrt{u_i(0)u_i(\frac\pi2)}}{||T_{\lambda_i}u_i||_{L^\infty}}.
$$
Furthermore $u_i(0)=u_i(\pi)$ implies
$w_i(\pi)=w_i(0)=w_i(\pi/2)$. Since  $\{w_i\}$ is a minimizing
sequence, we have
$$
\int_0^{2\pi}(w_i)^2_\theta-(w_i)^2d\theta\le 0
$$
for large $i$. It follows from $||w_i||_{L^\infty}=1$ that
$||w_i||_{H^1}\le C$. Hence $w_i\rightharpoonup w_0$ in $H^1$.
From  Sobolev embedding theorem we have that $w_i\to w_0$ in
$C^{0,\alpha}$ for any $\alpha\in(0,\frac12)$. If $\min w_0>0$,
then we know  that $\inf_{H^1_s}J(u)$ is bounded below  and  $w_0$
is a minimizer. Therefore we only need to consider the case when
$\min w_0=0$.

We  need  the following lemma to locate the zero points of $w_0$.
\begin{lemma}
\label{eqpi}
If $I=(a,b)$ is an open interval satifying $w_0(a)=w_0(b)=0$ and
$w_0(\theta)>0$ in $(a,b)$, then $b-a=\pi$.
\end{lemma}
\begin{proof}
If $b-a>\pi$, let $\delta=b-a-\pi$ and $J=(a+\frac\delta3,
b-\frac\delta3)$. Then $|J|=b-a-\frac{2\delta}{3}>\pi$. Therefore we
can choose $\beta>0$ such that $\cos(\theta+\beta)>0$ for
$\theta\in S^1\setminus J$, thus
$$
\int_{S^1\setminus J}\frac{\cos(\theta+\beta)}{w_i^3}d\theta\to
\infty \  \mbox{ as }i\to\infty,
$$
while $\int_J\frac{\cos(\theta+\beta)}{w_i^3}d\theta$ stays bounded.
So, as $i\to\infty$
$$
\int_0^{2\pi}\frac{\cos(\theta+\beta)}{w_i^3}d\theta\to \infty.
$$
This contradicts with
$$
\int_0^{2\pi}\frac{\cos(\theta+\beta)}{w_i^3}d\theta=0.
$$

On the other hand, we write the open set
$\{\theta:w_0(\theta)>0\}$ as a union of open intervals
$$
\{\theta:w_0(\theta)>0\}=\cup_i I_i.
$$
From the above argument we know that $|I_i|\le \pi$ for all $i$.
If one of these intervals, $I_{i_0}$ has length strictly less than
$\pi$, then
$$
\int_{I_{i_0}}((w_0')^2-w_0^2)d\theta>0,
$$
and
$$
\int_0^{2\pi}((w_0')^2-w_0^2)d\theta\ge
\int_{w_0>0}((w_0')^2-w_0^2)d\theta
=\sum_i\int_{I_i}((w_0')^2-w_0^2)d\theta>0.
$$
It follows that $J(w_i)>0$ for large $i$. This contradicts with the
fact that
$$
\inf_{H^1_s}J(u)\le J(1)=-4\pi^2.
$$
\end{proof}

Suppose $w_0(\theta_0)=\min w_0=0$ for some $\theta_0 \in [0,
\pi]$. Since $w_0(0)=w_0(\pi)=w_0(\pi/2) \ge w_0(3\pi/2)$, the
above lemma implies that if $w_0(0)=0$, then $w_0\equiv 0.$  This
contradicts with the fact of $\max w_0=1$. It follows that
$\theta_0 \ne 0$, $w_0(\theta_0)=w_0(\theta_0+\pi)=0$ and
$w_0(\theta)>0$ for $\theta\ne\theta_0, \theta_0+\pi$.

For each $i$, if $w_i(\theta_i)=\min w_i(\theta)$, then up to a
subsequence, we may assume, without loss of generality, that
$\theta_i\to \theta_0+\pi$. Furthermore, by rotating we may assume
that $\theta_0=0$, and $\theta_i=\theta_0+\pi=\pi$, that is, $w_i$
attains its minimum at $\pi$, $w_i\to w_0$ and $w_0>0$ for $\theta
\in [0, 2\pi]\setminus \{0, \pi\}.$

Let
$$
\tau_i=\sqrt{\frac{w_i(\frac\pi2)}{w_i(0)}},
\quad \sigma_i(\theta)=\sigma_{\tau_i}(\theta),
\mbox{  and  }
v_i=\frac{T_{\tau_i}w_i}{||T_{\tau_i}w_i||_{L^\infty}}.
$$
Again $\{v_i\}$ is a minimizing sequence and $||v_i||_{H^1}\le C$.
Therefore, up to a subsequence, $\{v_i\}$ converges weakly in
$H^1$ to some $v_0\in H^1$ and $v_i\to v_0$ in $C^{0,\alpha}$ for
any $\alpha\in (0,\frac12)$.

It is easy to see that
$$
v_i(0)=v_i(\frac\pi2), \,\lim_{i\to\infty}\frac{v_i(\frac{3\pi}2)}
{v_i(\frac\pi2)}=\frac{w_0(\frac{3\pi}{2})}{w_0(\frac\pi2)}\le 1,
\mbox{ and }v_i(\pi)\le v_i(0).
$$
Hence
\begin{equation}
\label{v0}
v_0(0)=v_0(\frac\pi2),\,\frac{v_0(\frac\pi2)}{v_0(\frac{3\pi}2)}\ge 1 \mbox{ and }v_0(\pi)\le v_0(0).
\end{equation}
Applying Lemma \ref{eqpi} to $v_0$, we see that $v_0(0)>0$,
otherwise $v_0(0)=v_0(\frac\pi2)=v_0(\pi)=v_0(\frac{3\pi}2)=0$,
which yields that $v_0\equiv 0$ by Lemma \ref{eqpi}, contradiction
to $\max v_0=1$. It follows from (\ref{v0}) that
$v_0(\frac\pi2)\ne 0$ , thus $v_0(\frac{3\pi}2)\ne0$. Applying
Lemma \ref{eqpi} again, we have $v_0(\pi)>0$. For any
$\theta\in[0,2\pi]\setminus \{\frac\pi2,\frac{3\pi}{2}\}$, since
$w_i$ attains its minimum at $\pi$, we obtain
\begin{align*}
v_i(\theta)=&\frac1{||T_{\tau_i}w_i||_{L^\infty}}w_i(\sigma_i(\theta))
\sqrt{\tau_i^2\cos^2\theta+\tau_i^{-2}\sin^2\theta}\\
\ge &\frac1{||T_{\tau_i}w_i||_{L^\infty}}w_i(\pi)\tau_i|\cos\theta|\\
=&v_i(\pi)|\cos\theta|\to v_0(\pi)|\cos\theta|>0
\end{align*}
Hence $v_0>0$, thus it  is a minimizer.

To complete the proof, we suppose that $v$ is a minimizer and observe
that
$$
H_s^1(S^1)= \{ u \in H^1(S^1) \ : \ u>0, \ \int_0^{2 \pi}
\cos \theta \cdot u^{-3} d \theta=\int_0^{2 \pi} \sin \theta \cdot
u^{-3} d \theta=0\}.
$$
Thus $v(\theta)$ satisfies the
Euler-Lagrange equation:
\begin{equation}\label{elab}
v_{\theta\theta}+v=\tau v^{-3}+\frac{a \cos \theta}{v^4}
+\frac{b\sin\theta}{v^4}=\tau v^{-3}+\frac{c \cos( \theta+\alpha)}{v^4},
\end{equation}
where $\tau,a, b$ are Lagrange multipliers, $c=\sqrt{a^2+b^2}$ and
$\gamma=\gamma(a,b)$. Multiplying both sides of (\ref{elab}) by
$\cos(\theta+\gamma)$ and integrating over $[0,2\pi]$, we have
$c=0$. Hence, $v$ satisfies
$$
v_{\theta\theta}+v=\tau v^{-3}.
$$
Due to Lemma \ref{eq}, we may assume $v(0)=v(\pi)$.
Define $w$ on $S^1$ by $w(\theta):=v(\frac\theta2)$,
$\theta\in[0,2\pi]$. One can check that $w\in H^1(S^1)$
satisfies
$$
w_{\theta\theta}+\frac14w=\frac14\tau w^{-2},
\qquad\theta\in S^1
$$
Using the same argument as in the proof of
Theorem \ref{thm1-1}, we obtain that
$$
w(\theta)=(\tau)^{1/4}\sqrt{\lambda^2\cos^2\frac{\theta-\alpha_1}2
+\lambda^{-2}\sin^2\frac{\theta-\alpha_1}2}.
$$
Therefore,
$$
v(\theta)=c\sqrt{\lambda^2\cos^2(\theta-\alpha)+\lambda^{-2}\sin^2(\theta-\alpha)}.
$$

\begin{remark}Our argument also yields another proof of the affine isoperimetric inequality.
In this regard, we refer \cite{CHLYZ} to interesting readers for the setting up.
\end{remark}

\section{Higher order derivatives}

We shall prove Theorem \ref{thm1-3} and \ref{thm1-4} in this section.

We first prove Theorem \ref{thm1-4}.
For any given $\lambda>0$ and $\alpha \in [0, 2\pi]$, let
$$
\Psi_{\lambda,\alpha}(\theta)=\left(
\lambda^2\cos^2\frac{\theta-\alpha}2
+\lambda^{-2}\sin^2\frac{\theta-\alpha}2\right)^{3/2},
$$

$$
\omega_{\lambda,\alpha}(\theta)=\alpha+\int_\alpha^\theta
{\Psi_{\lambda,\alpha}}^{-\frac23}d\theta
=
\left\{
\begin{array}{ll}
\alpha+2\arctan(\lambda^{-2}\tan\frac{\theta-\alpha}{2}), &\mbox{if} \ 0\le \theta -\alpha \le  {\pi}\\
\alpha+2\arctan(\lambda^{-2}\tan\frac{\theta-\alpha}{2})+2\pi, &\mbox{if} \ \pi<\theta -\alpha \le  {2 \pi}
\end{array}
\right.
$$
and
$$
({\mathcal T}_{\lambda ,\alpha}u)(\theta)
=u(\omega_{\lambda,\alpha}(\theta))\Psi_{\lambda,\alpha}
(\theta).
$$
From conformal invariant properties, we know that
$$
\int_0^{2\pi}(u^2_{\theta\theta}-\frac52u^2_\theta
+\frac9{16}u^2)d\theta=
\int_0^{2\pi}(({\mathcal T}_{\lambda, \alpha}u)^2_{\theta\theta}
-\frac52({\mathcal T}_{\lambda,\alpha} u)^2_\theta
+\frac9{16}({\mathcal T}_{\lambda,\alpha} u)^2)d\theta
$$
and
$$
\int_0^{2\pi}u^{-2/3}d\theta=\int_0^{2\pi}
({\mathcal T}_{\lambda,\alpha}u)^{-2/3}d\theta.
$$

\begin{lemma}
For any $u\in H^2(S^1)$, $u>0$, there exist $\alpha\in [0, 2\pi)$ and
$\lambda\ge 1$, such that
$$
\int_0^{2\pi}({\mathcal T}_{\lambda,\alpha}u)(\theta)\cos\theta
d\theta=\int_0^{2\pi}({\mathcal T}_{\lambda,\alpha} u)(\theta)
\sin\theta d\theta=0.
$$
\end{lemma}
\begin{proof}
For $\alpha\in S^1$, $\lambda\ge 1$, define
$$
I_{\lambda,\alpha}(u)=\left(
\int_0^{2\pi}({\mathcal T}_{\lambda,\alpha}u)(\theta)
\cos\theta d\theta,
\int_0^{2\pi}({\mathcal T}_{\lambda,\alpha}u)(\theta)
\sin\theta d\theta
\right)\in{\mathbf R}^2.
$$
Suppose $\exists$ $u>0$, such that $I_{\lambda,\alpha}\ne (0,0)$
for all $(\lambda, \alpha)\in [1,\infty)\times S^1$. Define
$G:[1,\infty)\times S^1 \to S^1$ by
$$
G(\lambda, \alpha)=\frac{I_{\lambda, \alpha}(u)}
{||I_{\lambda, \alpha}(u)||_{{\bf R}^2}},
$$
where $||{\bf X}||_{{\bf R}^2}$ is the length of vector ${\bf X} \in {\bf R}^2$.
One can easily check that $G$ is continuous.

When $\lambda=1$, $\omega_{1,\alpha}(\theta)=\theta$,
$\Psi_{1,\alpha}(\theta)=1$.
Therefore
\begin{align*}
\int_0^{2\pi}({\mathcal T}_{1,\alpha} u)(\theta)
\cos\theta d\theta
&=\int_0^{2\pi}u(\theta)\cos\theta d\theta,\\
\int_0^{2\pi}({\mathcal T}_{1,\alpha} u)(\theta) \sin\theta
d\theta &=\int_0^{2\pi}u(\theta)\sin\theta d\theta.
\end{align*}
That is $G(1,\cdot)$ maps $S^1$ to single point on $S^1$.

On the other hand, as $\lambda\to \infty$, $\omega(\theta)
\to \alpha$ pointwisely for $\theta\in[\alpha,\alpha+\pi)$ and
$\omega(\theta)\to \alpha+2\pi$ pointwisely for $\theta\in
(\alpha+\pi,\alpha+2\pi]$. Therefore  as $\lambda\to \infty$
\begin{align*}
\lambda^{-3}\int_0^{2\pi}({\mathcal T}_{\lambda,\alpha}u)
(\theta)\cos\theta d\theta
&=\lambda^{-3}\int_0^{2\pi}u(\omega_{\lambda,\alpha}(\theta))
\Psi_{\lambda,\alpha}(\theta)\cos\theta d\theta\\
&=\int_0^{2\pi}u(\omega_{\lambda,\alpha}(\theta))
(\cos^2\frac{\theta-\alpha}2+\frac1{\lambda^4}
\sin^2\frac{\theta-\alpha}2)^{\frac32}\cos\theta d\theta\\
&\to \frac85u(\alpha)\cos\alpha
\end{align*}
and
$$
\lambda^{-3}\int_0^{2\pi}({\mathcal T}_{\lambda,\alpha}u)
(\theta)\sin\theta d\theta\to
\frac85u(\alpha)\sin\alpha.
$$
Hence
$$
\lim_{\lambda\to\infty}G(\lambda,\alpha)=(\cos\alpha,\sin\alpha).
$$
That is $G(\infty,\cdot)=\mbox{Id}$, the identity map.
Since $\pi_1(S^1) \ne 0$, $\mbox{Id}$ is not homotopy to the constant map.
Contradiction.
\end{proof}
\begin{remark}
There are some similarities between the above lemma and other
rearrangement lemma (see, for example, Lemma 2 on Page 85 in
\cite{CA}).
\end{remark}

For any $u\in H^1(S^1)$, $u>0$, let
$$
F(u):=\int_0^{2\pi}(u_{\theta\theta}^2-\frac52u_\theta^2
+\frac9{16}u^2)d\theta\cdot\left(\int_0^{2\pi}u^{-2/3}
d\theta\right)^3.
$$
If  $\{u_i\}$ is a minimizing sequence of $F$ with
$\int_0^{2\pi}u_i^{-2/3}d\theta=1$, due to the above lemma, we may
assume, without loss of generality, that $u_i$ satisfies
$$
\int_0^{2\pi}u_i\cos\theta d\theta
=\int_0^{2\pi}u_i\sin\theta d\theta=0.
$$
Consider the Fourier expansion of $u_i$:
$$
u_i(\theta)=c_{i,0}+\sum_{k=2}^{\infty}(a_{i,k}\cos k\theta
+b_{i,k}\sin k\theta).
$$
It follows that
\begin{align}
\label{poF}
F(u_i)=&\pi\sum_{k=2}^\infty k^4(a_{i,k}^2+b_{i,k}^2)
-\frac{5\pi}2\sum_{k=2}^\infty k^2(a_{i,k}^2+b_{i,k}^2)
+\frac{9\pi}8 c_{i,0}^2
+\frac{9\pi}{16}\sum_{k=2}^\infty (a_{i,k}^2+b_{i,k}^2)\\
\ge & \frac{\pi}4\sum_{k=2}^\infty k^4(a_{i,k}^2+b_{i,k}^2)
+\frac{9\pi}{8}c_{i,0}^2 \ge C
\int_0^{2\pi}(((u_i)_{\theta\theta})^2+u_i^2)d\theta \nonumber
\end{align}
 for some positive constant $C$. Since
$\{u_i\}$ is a minimizing sequence with $\int_0^{2\pi}
u_i^{-2/3}d\theta=1$, $F(u_i)$ is bounded. It follows that
$||u_i||_{H^2}\le C$. Thus $u_i\rightharpoonup u_0$ weakly in
$H^2(S^1)$. Using Sobolev embedding theorem we obtain that
$u_0,u_i\in C^{1,\frac12}(S^1)$ and $u_i\to u_0$ in
$C^{1,\alpha}(S^1)$ for any $\alpha\in(0,\frac12)$.

If $u_0$ vanishes at some point, then
$\int_0^{2\pi}u_0^{-2/3}d\theta=\infty$ since
$u_0\in C^{1,\frac12}(S^1)$. On the other hand, Fatou's Lemma
implies
$$
\int_{0}^{2\pi}u_0^{-2/3}(\theta)d\theta\le
\lim_{i\to\infty}\int_{0}^{2\pi} u_i^{-2/3}(\theta)d\theta=1.
$$
Contradiction. Therefore $u_0>0$, thus it is a minimizer.

Suppose $u$ is a minimizer, then it satisfies the
Euler-Lagrange equation
$$
u_{\theta\theta\theta\theta}+\frac52u_{\theta\theta}
+\frac9{16}u=\tau u^{-\frac 53},
$$
where $\tau$ is a positive constant, since $\inf F(u)$ is
positive(by (\ref{poF})). Using the Green's function of the
operator $P_g$, we obtain that
\begin{equation}\label{green}
u(\theta)=c\int_0^{2\pi}\tau u^{-\frac 53}(\varphi)
|\sin\frac{\theta-\varphi}2|^3d\varphi.
\end{equation}
Define $v:{\mathbf R}\to{\mathbf R}$ by
\begin{equation}\label{pullback}
v(y):=u(\Phi^{-1}(y))\varphi^3(y),
\end{equation}
where $\Phi:S^1\to {\mathbf R}^1$ is
the stereographic projection:
$$
\Phi(\theta)=\tan\frac\theta 2.
$$
and
$$
\varphi(y):=\left(\frac{1+y^2}{2}\right)^{1/2}.
$$
It follows from (\ref{green}) that $v(y)$ satisfies
$$
v(y)=c\int_{-\infty}^\infty\tau v^{-\frac 53}|x-y|^3dx.
$$
Using Theorem 1.5 in Li \cite{Li} we obtain that
$$
v(y)=c\left(\frac{\lambda^2+\lambda^{-2} (y-a)^2}2\right)^{3/2},
$$
for some $c,\lambda>0$. Therefore
$$
u(\theta)=c\left(\lambda^2\cos^2\frac{\theta-\alpha}2
+\lambda^{-2}\sin^2\frac{\theta-\alpha}2\right)^{3/2}.
$$
This completes the proof of Theorem \ref{thm1-4}.

\medskip

\begin{remark}
\label{rem6-1} The same argument yields another proof of Theorem
\ref{thm1-2}.
\end{remark}
\begin{remark}
Theorem \ref{thm1-4} implies that $(S^1,g)$ has constant $Q$ curvature
if and only if it has constant $4$-scalar curvature.
\end{remark}

Finally, we shall turn to the  proof of Theorem \ref{thm1-3}. For any $\lambda>0$, let
$$
\Gamma_\lambda(\theta)=\left(\lambda^2\cos^2\theta
+\lambda^{-2}\sin^2\theta\right)^{3/2},
\quad \sigma_\lambda(\theta)=\int_0^\theta
\Gamma_\lambda^{-2/3}d\theta,
$$
and
$$
({\bf T}_\lambda u)(\theta)
:=u(\sigma_\lambda(\theta))\Gamma_\lambda(\theta).
$$

By conformal invariant property, we know that
$$
\int_0^{2\pi}(u^2_{\theta\theta}-10u^2_\theta
+9u^2)d\theta=
\int_0^{2\pi}(({\bf T}_\lambda u)^2_{\theta\theta}
-10({\bf T}_\lambda u)^2_\theta
+9({\bf T}_\lambda u)^2)d\theta
$$
and
$$
\int_0^{2\pi}u^{-2/3}d\theta=\int_0^{2\pi}
({\bf T}_\lambda u)^{-2/3}d\theta.
$$
Also one can easily check: for any $\alpha \in[0, 2\pi)$,
$$
\int_0^{2\pi}\frac{\cos^3(\theta+\alpha)}{({\bf T}_\lambda u)^{5/3}}d\theta
=(\lambda^{-2}\cos^2\alpha+\lambda^2\sin^2\alpha)^{3/2}
\int_0^{2\pi}\frac{\cos^3(\theta+\tilde{\alpha})}{u^{5/3}}d\theta,
$$
where $\tilde{\alpha}=\sigma^{-1}_\lambda(\alpha)$.
For any $u\in H^2(S^1)$, $u>0$, let
$$
{\bf F}(u):=\int_0^{2\pi}(u_{\theta\theta}^2-10u_\theta^2
+9u^2)d\theta\cdot\left(\int_0^{2\pi}u^{-2/3}
d\theta\right)^3,
$$
and
$$H_s^2(S^1):=\{ u \in H^2(S^1) \ : \ \int_0^{2\pi}\frac{\cos^3(\theta+\alpha)}{ u^{5/3}}d\theta=0
\ \mbox{for} \ \mbox{all} \ \alpha\in [0, 2\pi)\}.$$ We first
prove the existence of a minimizer of $\inf_{H_s^2(S^1)} {\bf
F}(u)$.

Suppose that $\{u_i\}$ is a minimizing sequence of
$\inf_{H_s^2(S^1)} {\bf F}(u)$, with $||u_i||_{L^{\infty}}=1$. The
following lemma implies $||u_i||_{H^2}$ is bounded, therefore
$u_i\rightharpoonup u_0$ in $H^2(S^1)$. Using Sobolev embedding
theorem we obtain that $u_i, u_0\in C^{1,\frac12}$, and $u_i\to
u_0$ in $C^{1,\alpha}$, for any $\alpha\in(0,\frac12)$.
\begin{lemma}\label{bound}
Suppose $\{v_i\}$ is a minimizing sequence of $\inf_{H_s^2(S^1)}
\bf F$ with $||v_i||_{L^\infty}=1$. Then $||v_i||_{H^2}$ is
bounded.
\end{lemma}
\begin{proof} We argue by contradiction.
Suppose that $||v_i||_{H^2}$ is unbounded. Up to a subsequence, we may assume
that $\lim_{i\to\infty}||v_i||_{H^2}=\infty$. Let
$w_i=v_i/||v_i||_{H^2}$. Then $\{w_i\}$ is also a minimizing sequence,
$w_i\rightharpoonup 0$ in $H^2(S^1)$ and  $||w_i||_{H^2}=1$.
It follows that  $w_i\to 0$
in $C^{1,\alpha}$ for any $\alpha\in(0,\frac12)$, thus  $\lim_{i\to\infty}\int_0^{2\pi}
(w_i^2+(w_i)_\theta^2)d\theta=0$. On the other hand, we know from $w_i \rightharpoonup 0$ in $H^2(S^2)$ that
$$
\lim_{i\to\infty}\int_0^{2\pi}(w_i)^{-2/3}d\theta=\infty.
$$
Note that
$$
\lim_{i\to\infty}\int_0^{2\pi}((w_i)_{\theta\theta}^2-10(w_i)_\theta^2
+9w_i^2)d\theta\ge 0,
$$
and $\{w_i\}$ is a minimizing sequence. We thus have
$$
\lim_{i\to\infty}\int_0^{2\pi}((w_i)_{\theta\theta}^2-10(w_i)_\theta^2
+9w_i^2)d\theta= 0.
$$
Therefore $\lim_{i\to\infty}\int_0^{2\pi}((w_i)_{\theta\theta}^2d\theta
=0$. Hence $\lim_{i\to\infty}||w_i||_{H^2}=0$. This contradicts with
$||w_i||_{H^2}=1$.
\end{proof}

If $u_0>0$, we know as in the proof of Theorem \ref{thm1-1} that
$u_0$ is a minimizer of $\inf_{H_s^2(S^1)} {\bf F}(u)$. The subtle
case is that $\min u_0=0$.

We need one more lemma.

\begin{lemma}\label{eqpi3}
If $I=(a,b)$ is an open interval satifying $u_0(a)=u_0(b)=0$ and
$u_0(\theta)>0$ in $(a,b)$, then $b-a=\pi$.
\end{lemma}
\begin{proof}
If $b-a>\pi$, let $\delta=b-a-\pi$ and $J=(a+\frac\delta3,
b-\frac\delta3)$. Then $|J|=b-a-\frac{2\delta}{3}>\pi$. Therefore we
may choose $\beta>0$ such that $\cos(\theta+\beta)>0$ for
$\theta\in S^1\setminus J$. Since $u_0\in C^{1,\frac12}$,
$$
\int_{S^1\setminus J}\frac{\cos^3(\theta+\beta)}{u_i^{5/3}}d\theta\to
\infty,\mbox{ as }i\to\infty,
$$
while $\int_J\cos^3(\theta+\beta)/u_i^{5/3}d\theta$ is
bounded. Therefore, as $i\to\infty$
$$
\int_0^{2\pi}\frac{\cos^3(\theta+\beta)}{u_i^{5/3}}d\theta\to \infty.
$$
Contradiction with the fact that $u_i\in H_s^2(S^1).$

On the other hand, write the open set $\{\theta:u_0(\theta)>0\}$ as
a union of open intervals:
$$
\{\theta:w_0(\theta)>0\}=\cup_i I_i.
$$
From the above argument we know that $|I_i|\le \pi$ for all $i$.
We claim that for interval $I_k$ with $|I_k|\le\pi$
$$
\int_{I_k}((u_0)_{\theta\theta}^2-10(u_0)_\theta^2+9u_0^2)
d\theta\ge0,
$$
and for $I_l$ with $|I_l|<\pi$,
$$
\int_{I_l}((u_0)_{\theta\theta}^2-10(u_0)_\theta^2+9u_0^2)
d\theta>0.
$$
Therefore, if for some $i_0$, $|I_{i_0}|< \pi$, we have
\begin{align*}
\int_0^{2\pi}((u_0)_{\theta\theta}^2-10(u_0)_\theta^2+9u_0^2)
d\theta\ge& \int_{u_0>0}((u_0)_{\theta\theta}^2-10(u_0)_\theta^2
+9u_0^2)d\theta\\
=&\sum_i\int_{I_i}((u_0)_{\theta\theta}^2-10(u_0)_\theta^2
+9u_0^2)d\theta\\
\ge&\int_{I_0}((u_0)_{\theta\theta}^2-10(u_0)_\theta^2
+9u_0^2)d\theta=C>0.
\end{align*}
Since $u_0\in C^{1,1/2}$, thus $\int_0^{2 \pi} u_i^{-2/3} d \theta \to \infty$.
It follows that
$$
\lim_{i\to\infty} {\bf F}(u_i)=\infty.
$$
Contradiction.

We return to the above claim. Suppose $|I_k|=\pi$. Without loss of generality
we may assume that $I_k=(0,\pi)$, i.e. $u_0(0)=0$, $u_0(\pi)=0$ and
$u_0(\theta)>0$ when $\theta\in(0,\pi)$. Since $u\ge 0$ and
$u\in C^{1,\frac12}$, we obtain that $u'(0)=u'(\pi)=0$. Therefore
if $v(\theta)=u(\frac\theta2)$ for $\theta\in[0,2\pi]$, then $v \in H^2(S^1)$.
It follows from Theorem \ref{thm1-4} that
$$
\int_{I_k}((u_0)_{\theta\theta}^2-10(u_0)_\theta^2+9u_0^2)d\theta
=16\int_0^{2\pi}(v_{\theta\theta}^2-\frac52v_\theta^2+\frac9{16}v^2)
d\theta\ge 0.
$$
For $|I_l|< \pi$, we assume that $I_l=(-a,a)$ with $a<\frac\pi2$. Let
$v(\theta)=u_0(\frac\theta2)$ and $w(y)=v(2\arctan y)
((1+y^2)/2)^{3/2}$(see (\ref{pullback})). Since
$u_0(a)=u_0(-a)=0$ and $u_0'(-a)=u_0'(a)=0$, we know that
$w(\pm \tan a)=w'(\pm \tan a)=0$ and $w(y)>0$ for $y\in(-\tan a,\tan a)$.
Hence
\begin{align*}
\int_{I_l}((u_0)_{\theta\theta}^2-10(u_0)_\theta^2+9u_0^2)d\theta
&=16\int_{-2a}^{2a}(v_{\theta\theta}^2-\frac52v_\theta^2
+\frac9{16}v^2)d\theta\\
&\ge \int_{-\tan a}^{\tan a}(w''(y))^2dy>0.
\end{align*}
\end{proof}

The rest of the proof of the existence of a minimizer will be similar
to the proof of the existence part of Theorem \ref{thm1-1}.

Using Lemma \ref{eq}, we may assume that $u_i(0)=u_i(\pi)$ and
$u_i(\pi/2)\ge u_i(3\pi/2)$. Let
$$
\lambda_i=\sqrt{\frac{u_i(\frac\pi2)}{u_i(0)}},
\mbox{ and }
w_i(\theta)=\frac{{\bf T}_{\lambda_i}u_i}{||{\bf T}_{\lambda_i}u_i
||_{L^\infty}}.
$$
It is easy to see that $\{w_i\}$ is also a minimizing sequence and
$$
w_i(0)=w_i(\frac\pi2)=\frac{\sqrt{u_i(0)u_i(\frac\pi2)}}
{||{\bf T}_{\lambda_i}u_i ||_{L^\infty}}.
$$
Furthermore $u_i(0)=u_i(\pi)$ imples $w_i(\pi)=w_i(0)=w_i(\pi/2)$.
Since  $\{w_i\}$ is a minimizing sequence with
$||w_i||_{L^\infty}=1$, Lemma \ref{bound} implies that
$||w_i||_{H^2}$ is bounded. Hence $w_i\rightharpoonup w_0$ in
$H^2$. It follows from  Sobolev embedding theorem that $w_i\to
w_0$ in $C^{1,\alpha}$ for any $\alpha\in(0,\frac12)$. If $\min
w_0>0$, then $w_0$ is a minimizer as being pointed out above.
Therefore we only need to consider the case of $\min w_0=0$.

Suppose $w_0(\theta_0)=\min w_0=0$ for some $\theta_0\in[0,\pi]$.
Since $w_0(0)=w_0(\pi)=w_0(\pi/2) \ge w_0(3\pi/2)$, Lemma
\ref{eqpi3} implies that if $w_0(0)=0$ then $w_0\equiv 0$, which
contradicts with the fact that $\max w_0=1$. Hence $\theta_0\ne
0$, $w_0(\theta_0) =w_0(\theta_0+\pi)=0$ and $w_0(\theta)>0$ for
$\theta\ne\theta_0, \theta_0+\pi$.

For each $i$, suppose $w_i(\theta_i)=\min w_i(\theta)$.
After rotating and up to a subsequence of $w_i$, we may assume that
$\theta_0=0$, and $\theta_i=\theta_0+\pi=\pi$, i.e.
$w_i$ attains its minimum at $\pi$, $w_i\to w_0$ and $w_0>0$ except at
$0$ and $\pi$.

Let
$$
\tau_i=\sqrt{\frac{w_i(\frac\pi2)}{w_i(0)}},
\quad \sigma_i(\theta)=\sigma_{\tau_i}(\theta),
\mbox{  and  }
v_i=\frac{{\bf T}_{\tau_i}w_i}{||{\bf T}_{\tau_i}w_i||_{L^\infty}}.
$$
Again $\{v_i\}$ is a minimizing sequence and $||v_i||_{H^2}\le C$
by Lemma \ref{bound}. Therefore, up to a subsequence, $\{v_i\}$
converges weakly in $H^2$ to some $v_0\in H^2$ and $v_i\to v_0$ in
$C^{1,\alpha}$ for any $\alpha\in (0,\frac12)$.

It is easy to see that
$$
v_i(0)=v_i(\frac\pi2), \,\lim_{i\to\infty}\frac{v_i(\frac{3\pi}2)}
{v_i(\frac\pi2)}=\frac{w_0(\frac{3\pi}{2})}{w_0(\frac\pi2)}\le 1,
\mbox{ and }v_i(\pi)\le v_i(0).
$$
Hence
\begin{equation}
\label{v03}
v_0(0)=v_0(\frac\pi2),\,\frac{v_0(\frac\pi2)}{v_0(\frac{3\pi}2)}
\ge 1\mbox{ and }v_0(\pi)\le v_0(0).
\end{equation}
Applying Lemma \ref{eqpi3} to $v_0$, we see that $v_0(0)>0$,
otherwize $v_0(0)=v_0(\frac\pi2)=v_0(\pi)=v_0(\frac{3\pi}2)=0$,
contradiction. It follows from (\ref{v03}) that $v_0(\frac\pi2)\ne
0$, thus $v_0(\frac{3\pi}2)\ne0$. Applying Lemma \ref{eqpi3}
again, we have $v_0(\pi)>0$. For any $\theta\in[0,2\pi]\setminus
\{\frac\pi2,\frac{3\pi}{2}\}$, since $w_i$ attains its minimum at
$\pi$, we obtain
\begin{align*}
v_i(\theta)=&\frac1{||{\bf T}_{\tau_i}w_i||_{L^\infty}}w_i(\sigma_i(\theta))
(\tau_i^2\cos^2\theta+\tau^{-2}\sin^2\theta)^{3/2}\\
\ge &\frac1{||{\bf T}_{\tau_i}w_i||_{L^\infty}}w_i(\pi)\tau_i|\cos\theta|^3\\
=&v_i(\pi)|\cos\theta|^3\to v_0(\pi)|\cos\theta|^3>0
\end{align*}
Hence $v_0>0$, thus it is a minimizer.

To classify the extremal function $v(\theta)$,  we again observe
that
\begin{align*}H_s^2(S^1)= \{ u \in H^2(S^1) \ :& \ u>0, \ \int_0^{2 \pi}
\cos 3\theta \cdot u^{-5/3} d \theta
=\int_0^{2 \pi} \sin 3\theta \cdot u^{-5/3} d \theta=0\\
&\mbox{and }\int_0^{2
\pi} \cos \theta \cdot u^{-5/3} d \theta=\int_0^{2 \pi} \sin
\theta \cdot u^{-5/3} d \theta=0\}.
\end{align*}

Thus $v(\theta)$
satisfies the Euler-Lagrange equation:
\begin{align}
\label{elab-1}
v_{\theta\theta \theta \theta}+10 v_{\theta\theta} +9v=&\tau v^{-\frac 53}
+\frac{a\cos3\theta}{v^{\frac83}}+\frac{b\sin3\theta}{v^{\frac 83}}
+\frac{c\cos\theta}{v^{\frac 83}}+\frac{d\sin \theta}{v^{\frac 83}}\\
=&\tau v^{-\frac 53}
+\frac{A\cos(3\theta+\alpha_1)}{v^{\frac83}}
+\frac{B\cos(\theta+\alpha_2)}{v^{\frac 83}}\nonumber
\end{align}
where $\tau,a, b, c, d$ are Lagrange multipliers and
$A=\sqrt{a^2+b^2}$, $B=\sqrt{c^2+d^2}$. We then follow the
argument given in \cite{CHLYZ}. Multiplying both sides of
(\ref{elab-1}) by $\cos(3\theta+\alpha_1)$ and
$\cos(\theta+\alpha_2)$ respectively and integrating over
$[0,2\pi]$, we obtain that
\begin{equation}\label{lsab-1}
\begin{cases}
&A\int_0^{2\pi}\displaystyle\frac{\cos^2(3\theta+\alpha_1)}{v^{8/3}}d\theta
+B\int_0^{2\pi}\displaystyle\frac{\cos(3\theta+\alpha_1)\cos(\theta+\alpha_2)}
{v^{8/3}}d\theta=0\\
&A\int_0^{2\pi}\displaystyle\frac{\cos(3\theta+\alpha_1)\cos(\theta+\alpha_2)}
{v^{8/3}}
d\theta +B\int_0^{2\pi}\displaystyle\frac{\cos^2(\theta+\alpha_2)}{v^{8/3}}
d\theta=0.
\end{cases}
\end{equation}
It follows from the H\"older inequality that the determinant of
the coefficient matrix of linear system (\ref{lsab-1}) is always
positive. Therefore (\ref{lsab-1}) has only the trivial solution
$A=B=0$. Hence, $v$ satisfies
$$
v_{\theta\theta \theta \theta}+10 v_{\theta\theta } +9v=\tau v^{-3},
$$
which implies that the symmetric $Q$-curvature of $v^{-\frac
43}g_s$ is a constant.  We hereby completes the proof of Theorem
\ref{thm1-3}.

\begin{remark}\label{rem6-2}
It can be  checked from the classification of constant {\it
1-scalar curvature} metrics that the symmetric $Q$-curvature of a
constant {\it 1-scalar curvature} metric is a constant.
Conversely, it is not clear whether it is true, even though we
suspect that it is the case. Technically, we do have difficulty to
classify the metrics with constant symmetric $Q$ curvature.  We
conjecture that if $u$ is a minimizer of $\inf_{H_s^2(S^1)} {\bf
F}(u)$, then
$$
u(\theta)=c\left(\lambda^2\cos^2(\theta-\alpha)
+\lambda^{-2}\sin^2(\theta-\alpha)\right)^{3/2},
$$
for some $\lambda,c>0$ and $\alpha \in [0, 2\pi)$, and the infimum
is ${\bf F}(u)=144\pi^4$.
\end{remark}
\begin{remark}
Direct computation shows that
$$
Q^\alpha_g=\frac{\alpha}3\Delta_g R^\alpha_g + (R^\alpha_g)^2,
$$
which implies
$$
\int Q^\alpha_g d S_g= \int (R^\alpha_g)^2 d S_g.
$$
Therefore using Theorem \ref{thm1-4} we obtain that with length
fixed, the infimum of $\int (k_g)^2 d S_g$ is attained by metric
with constant $4$-scalar curvature. In addition, it can be shown
that $\int Q_g d S_g$, therefore $\int (k_g)^2 d S_g$ is
increasing and bounded above along the $Q$-curvature flow
(\ref{*q}).
\end{remark}

\end{document}